\documentclass{article}
\usepackage{amsmath, epsfig,amssymb, multicol}
\usepackage[active]{srcltx}
\newtheorem{lemma}{Lemma}
\newtheorem{theorem}{Theorem}

\newtheorem{corollary}{Corollary}

\newtheorem{remark}{Remark}

\renewcommand{\L}{{\cal L}}

\renewcommand{\S}{{\cal S}}
\newcommand{\T}{{\cal T}}
\renewcommand{\P}{{\cal P}}

\renewcommand{\d}{{\bf d}}
\newcommand{\1}{{\bf 1}}
\newcommand{\tr}{{\rm Trace}}
\newcommand{\G}{{\cal G}}

\newcommand{\E}{{\rm E}}
\newcommand{\Var}{{\rm Var}}
\newcommand{\Cov}{{\rm Cov}}

\newcommand{\Z}{{\mathbb Z}}

\renewcommand{\Pr}{{\rm Pr}}

\newcommand{\vol}{{\rm vol}}

\title{Loose  Laplacian spectra of  random hypergraphs }
\author{Linyuan Lu
\thanks{University of South Carolina, Columbia, SC 29208,
({\tt lu@math.sc.edu}). This author was supported in part by NSF
grant DMS 1000475. }
  \and Xing Peng
\thanks{University of South Carolina, Columbia, SC 29208,
({\tt pengx@mailbox.sc.edu}).This author was supported in part by
NSF grant  DMS 1000475. }}
\begin{document}
\maketitle

\begin{abstract}
Let $H=(V,E)$ be an $r$-uniform hypergraph with the vertex set $V$ and the
edge set $E$. For $1\leq s \leq r/2$, we define a weighted graph $G^{(s)}$
on the vertex set ${V\choose s}$ as follows. Every pair of $s$-sets $I$ and $J$
is associated with a weight $w(I,J)$, which is the number of edges in $H$ passing
through $I$ and $J$ if $I\cap J=\emptyset$, and $0$ if $I\cap J\not=\emptyset$.
The $s$-th Laplacian $\L^{(s)}$ of $H$ is defined to be the  normalized Laplacian of $G^{(s)}$.
The eigenvalues of $\mathcal L^{(s)}$
are listed as $\lambda^{(s)}_0, \lambda^{(s)}_1, \ldots,
\lambda^{(s)}_{{n\choose s}-1}$ in non-decreasing order.
Let $\bar\lambda^{(s)}(H)=\max_{i\not=0}\{|1-\lambda^{(s)}_i|\}$.
The parameters $\bar\lambda^{(s)}(H)$ and $\lambda^{(s)}_1(H)$,
which were introduced in our previous paper, have
a number of connections to the mixing rate of high-ordered random walks, the
generalized distances/diameters, and the edge expansions.

For $0< p<1$,  let $H^r(n,p)$ be a random $r$-uniform hypergraph over
$[n]:=\{1,2,\ldots, n\}$,  where each $r$-set of $[n]$
 has probability $p$ to be an edge independently.
For $1 \leq s \leq r/2$, $p(1-p)\gg \frac{\log^4 n}{n^{r-s}}$, and $1-p\gg
\frac{\log n}{n^2}$,  we prove that
 almost surely
$$\bar\lambda^{(s)}(H^r(n,p))\leq
 \frac{s}{n-s}+
\left(3+o(1)\right)\sqrt{\frac{1-p}{{n-s\choose r-s}p}}.$$ We also
prove that the empirical distribution of the eigenvalues of
$\L^{(s)}$ for $H^r(n,p)$ follows the Semicircle Law if
$p(1-p)\gg \frac{\log^{1/3} n}{n^{r-s}}$ and $1-p\gg
\frac{\log n}{n^{2+2r-2s}}$. 
\end{abstract}

\section{Introduction}
The spectrum of the adjacency matrix (and/or the Laplacian matrix) of
a random graph was well-studied in the literature \cite{akv, flv1, flv2,
  og,o, dj,fo,fks,fk}.  Given a graph $G$, let $\mu_1(G),\ldots,
\mu_n(G)$ be the eigenvalues of the adjacency matrix of $G$ in the
non-decreasing order, and $\lambda_0(G), \ldots,$ $\lambda_{n-1}(G)$
be the eigenvalues of (normalized) Laplacian matrix of $G$
respectively.  Let $G(n,p)$ be the Ed\H{o}s-R\'enyi random graph
model.  F\"uredi and Koml\'os \cite{fk} showed that if $np(1-p)\gg
\log^6 n$ then almost surely $\mu_n=(1+o(1))np$ and
$\max\{-\mu_1,\mu_{n-1}\}\leq (2+o(1))\sqrt{np(1-p)}$. The results
are extended to sparse random graphs \cite{fo, ks} and general
random matrices \cite{dj, fk}. Alon, Krivelevich, and Vu \cite{akv}
proved the concentration of the $s$-th largest eigenvalue of a
random symmetric  matrix with independent random entries of absolute
value at most $1$.
 Friedman (in a
series of papers \cite{fks, friedman, friedman2}) proved that the
second largest eigenvalue of random $d$-regular graphs is almost
surely $(2+o(1))\sqrt{d-1}$ for any $d\geq 4$.
 Chung, Lu, and Vu \cite{flv2} studied the Laplacian
eigenvalues of random graphs with given expected degrees; their
results were supplemented by Coja-Oghlan \cite{og,o} for much sparser
random graphs.

In this paper, we study the spectra of the Laplacians of random
hypergraphs.  Laplacians for regular hypergraphs were first
introduced by Chung \cite{fan2}  using the homology approach.
Rodr\'iguez \cite{rod1, rod2} treated a hypergraph as a multi-edge
graph and then defined its Laplacian to be the Laplacian of the
corresponding multi-edge graph. Inspired by these work, we
\cite{rwalk} introduced  the generalized Laplacian eigenvalues of
hypergraphs through high-ordered random walks. Let $H=(V,E)$ be an
$r$-uniform hypergraph on $n$ vertices.  We can associate $r-1$
Laplacians $\L^{(s)}$ ($1\leq s \leq r-1$) to $H$; roughly speaking,
$\L^{(s)}$ captures the incidence relations between $s$-sets and
edges in $H$. Our definition of the Laplacian at the spacial case
$s=1$ is the same as the Laplacian considered by Rodr\'iguez
\cite{rod1, rod2}.  The $s$-th Laplacian is {\em loose} if $1 \leq s
\leq r/2$, and is {\em tight} if $r/2<s\leq r-1$.  Here we consider
only the spectra of loose Laplacians.

For $1\leq s \leq r/2$, we consider an auxiliary weighted graph
$G^{(s)}$ defined as follows: the vertex set of $G^{(s)}$ is
${V\choose s}$ while the weight function $W\colon {V\choose
  s}\times {V\choose s} \to \Z$ is defined as

\begin{equation}
  \label{eq:1}
W(S,T)=\left \{
  \begin{array}[c]{ll}
\left|\{F\in E(H)\colon  S\cup T\subset F\}\right| &
\mbox{ if } S\cap T=\emptyset;\\
0 & \mbox{ otherwise.}
  \end{array}
\right.
\end{equation}

The {\em $s$-th Laplacian} of $H$, denoted by $\L^{(s)}$, is the
normalized Laplacian of $G^{(s)}$.  For any $s$-set $S$, let $d_S$ be
the number of edges in $H$ passing through $S$; the degree of $S$ in
$G^{(s)}$ is ${r-s\choose s}d_S$. Let $D$ be the diagonal matrix of
the degrees $\{d_S\}$ and $W$ be the weight
matrix $\{w(S,T)\}$.  Note that $T:={r-s\choose s}D$ is the diagonal matrix of degrees
in $G^{(s)}$. We have
\begin{equation}
  \label{eq:2}
  \L^{(s)}= I - T^{-1/2}W  T^{-1/2}.
\end{equation}

The eigenvalues of $\mathcal L^{(s)}$
are listed as $\lambda^{(s)}_0, \lambda^{(s)}_1, \ldots,
\lambda^{(s)}_{{n\choose s}-1}$ in non-decreasing order.
We have
\begin{equation}
\label{eq:3}
0=\lambda^{(s)}_0\leq  \lambda^{(s)}_1\leq \cdots \leq
\lambda^{(s)}_{{n\choose s}-1}\leq 2.
\end{equation}
The first non-trivial eigenvalue $\lambda^{(s)}_1>0$ if and only if
 $G^{(s)}$ is connected. When this occurs, we say
$H$ is {\it $s$-connected}. The diameter of $G^{(s)}$ is called
the {\em $s$-th diameter} of $H$.
The largest eigenvalue $\lambda^{(s)}_{{n\choose s}-1}$
is also denoted by $\lambda^{(s)}_{max}$. The (Laplacian) spectral radius,
denoted by $\bar \lambda^{(s)}$, is the maximum of $1-\lambda^{(s)}_1$
and $\lambda^{(s)}_{max}-1$.

This definition differs slightly with the one in \cite{rwalk}, where the vertex set of
the auxiliary graph (denoted by $G^{(s)'}$) is the set of all distinct
$s$-tuples instead. Note that $G^{(s)'}$ is the blow-up of
$G^{(s)}$. Their Laplacian spectra differ only by the multiplicity of
$1$'s. Therefore, two different definitions give the same values  of
 $\lambda^{(s)}_1$, $\lambda^{(s)}_{max}$, and $\bar\lambda^{(s)}$.

For different $s$, the following inequalities were proved in \cite{rwalk}.
\begin{eqnarray}
\label{labmda1} \lambda_1^{(1)} \geq \lambda_1^{(2)} \geq \ldots
\geq
\lambda_1^{(\lfloor r/2 \rfloor)};\\
\label{lambdamax} \lambda_{\max}^{(1)} \leq \lambda_{\max}^{(2)}
\leq \ldots \leq \lambda_{\max}^{(\lfloor r/2 \rfloor)}.
\end{eqnarray}

The $s$-th Laplacian has a number of connections to the mixing rate of
high-ordered random walks, the generalized distances/diameters, and
the edge expansions. Here we list some applications, which are
similar to results in \cite{rwalk}, and results for graphs \cite{fan1, fan3, fan4,fan5,
fan6, fan7}.
\begin{description}
\item[Random $s$-Walks:] The mixing rate of the random $s$-walk on $H$ is at most
$\bar\lambda^{(s)}$.
\item[The $s$-Diameter:] The $s$-diameter of $H$ is at most
$$
\left\lceil \frac{\log \frac{|E(H)|{r\choose s}}{\delta}}
  {\log
  \frac{\lambda_{\max}^{(s)}+\lambda_1^{(s)}}{\lambda_{\max}^{(s)}-\lambda_1^{(s)}}}\right\rceil.
$$
Here $\delta=\min_{S\in {V\choose s}} d_S$ is the minimum degree among all $s$-sets.
\item [Edge expansion:] For $1\leq t\leq s \leq \frac{r}{2}$, $\S\subset {V\choose t}$, and
$\T\subset {V\choose t}$, define
$$E(\S,\T)=\{F\in E(H)\colon \exists S\in \S,\exists T\in \T
\mbox{ such that } S \cap T=\emptyset,\mbox{ and }
S\cup T\subset F\},$$
$$e(\S,\T)=\frac{|E(\S,\T)|}{\left|E({V\choose s}, {V\choose t})\right|},
$$
$$e(\S)=\frac{\sum_{S\in \S}d_S}{\sum_{S\in {V\choose s}}d_S},$$
$$e(\T)=\frac{\sum_{T\in \T}d_T}{\sum_{T\in {V\choose t}}d_T}.$$
Then we have
$$|e(\S,\T)-e(\S)e(\T)|\leq
\bar \lambda^{(s)}\sqrt{e(\S)e(\T)e(\bar \S) e(\bar \T)}.$$
\end{description}
The proofs of these claims are very similar to those in \cite{rwalk} and are omitted here.

Our first result is the eigenvalues of the $s$-th Laplacian of the complete $r$-uniform hypergraph $K^r_n$.
\begin{theorem}\label{t0}
  Let $K^r_n$ be the complete $r$-uniform hypergraph on $n$ vertices.
For $1 \leq s \leq r/2$,
the eigenvalues of $s$-th Laplacian of $K^r_n$ are given by
$$
1- \frac{(-1)^{i}\binom{n-s-i)}{s-i}}{{n-s \choose s}}
\mbox{ with multiplicity } {n\choose i}-{n\choose i-1}
\mbox{ for }
0\leq i \leq s. $$
\end{theorem}
Here we point out an application of this theorem to the
celebrated Erd\H{o}s-Ko-Rado Theorem, which states ``if the $n\geq 2s$, then
the size of the maximum intersecting family of $s$-sets
in $[n]$ is at most ${n-1 \choose s-1}$.'' (The theorem was originally proved by
Erd\H{o}s-Ko-Rado \cite{EKR} for sufficiently large $n$; the simplest proof was due to
 Katona \cite{katona}.) Here we present a proof adapted from Calderbank-Frankl \cite{cf},
where they use the eigenvalues of Kneser graph instead. (The relation between
$\L^{(s)}(K^r_n)$ and the Laplacian of the Kneser graph is explained in section 2.)

It suffices to show for any intersecting family $U$ of $s$-sets, $|U|\leq {n-1\choose s-1}$.
 Note that $U$ is an independent set
of $G^{(s)}(K^r_n)$. Restricting to $U$,
 $\L^{(s)}(K^r_n)$ becomes an identity matrix; whose eigenvalues are all equal to $1$.
By Cauchy's interlace theorem, we have
\begin{equation}
  \label{eq:il}
\lambda_k^{(s)}\leq 1 \leq \lambda_{{n\choose s}-|U|+k}^{(s)}
\end{equation}
for $0\leq k\leq |U|-1$.
Let $N^+$ (or $N^-$) be the number
of eigenvalues of $\L^{(s)}(K^r_n)$ which is $\geq 1$ (or $\leq 1$) respectively.
Inequality (\ref{eq:il}) implies that $|U|\leq N^+$ and $|U|\leq N^-$.
By Theorem \ref{t0}, $N^+=\sum_{i=0}^{\lfloor (s-1)/2\rfloor}
\left( {n\choose 2i+1}- {n\choose 2i}\right)$
and $N^-=\sum_{i=0}^{\lfloor s/2\rfloor}
\left( {n\choose 2i}- {n\choose 2i-1}\right)$. We have
$$|U|\leq \min\{N^+,N^-\}=\sum_{i=0}^{s-1} (-1)^{s-1-i}{n\choose i}={n-1\choose s-1}.$$

For $0< p <1$,  let $H^r(n,p)$ be a random $r$-uniform hypergraph over
$[n]=\{1,2,\ldots, n\}$,  where each $r$-set of $[n]$
 has probability $p$ to be an edge independently.
We can estimate the Laplacian spectrum of $H^r(n,p)$ using
the Laplacian spectrum of $K^r_n$ as follows.

 \begin{theorem} \label{t1}
Let $H^r(n,p)$ be a random $r$-uniform hypergraph.
For $1 \leq s \leq r/2$, if $p(1-p)\gg \frac{\log^4 n }{n^{r-s}}$ and
$1-p\gg \frac{\log n}{n^2}$, then  almost surely
the $s$-th spectral radius $\bar\lambda^{(s)}(H^r(n,p))$
satisfies
\begin{equation}
  \label{eq:barlambda}
\bar\lambda^{(s)}(H^r(n,p))\leq \frac{s}{n-s}+
\left(3+o(1)\right)\sqrt{\frac{1-p}{{n-s\choose r-s}p}}.
\end{equation}
Moreover, for $1\leq k \leq {n\choose s}-1$,  almost surely we have
\begin{equation}
  \label{eq:diffub}
|\lambda_k^{(s)}(H^r(n,p))-\lambda_k^{(s)}(K^r_n)|\leq \left(3+o(1)\right)\sqrt{\frac{1-p}{{n-s\choose r-s}p}}.
\end{equation}
 \end{theorem}

Note that $G(n,p)$ is a special case of $H^r(n,p)$ with $r=2$.
By choosing $s=1$, Theorem \ref{t1} implies that
\begin{equation}
  \label{eq:np1}
\bar\lambda(G(n,p))\leq (3+o(1))\sqrt{\frac{1-p}{(n-1)p}}
\quad\mbox{ for } p(1-p)\gg \frac{\log^4 n}{n}.
\end{equation}

  Chung, Lu, and Vu's result\cite{flv2},
when restricted to $G(n,p)$, implies
\begin{equation}
  \label{eq:np2}
\bar\lambda(G(n,p))\leq (4+o(1))\frac{1}{\sqrt{np}}  \quad\mbox{ for } 1-\epsilon \geq p\gg \frac{\log^6 n}{n}.
\end{equation}

Inequality (\ref{eq:np1}) has
a smaller constant and works for a larger range of $p$ than inequality (\ref{eq:np2}).

F\"uredi and Koml\'os \cite{fk} proved the empirical distribution of  the
eigenvalues of $G(n,p)$ follows the Semicircle Law. Chung, Lu, and Vu \cite{flv2} proved
a similar result for the random graphs with  given expected degrees. Here
we prove a similar result for random hypergraphs.

 \begin{theorem} \label{t2}
 For $1\leq s \leq r/2$,   if $p(1-p)\gg \frac{\log^{1/3} n}{n^{r-s}}$ and $1-p\gg
\frac{\log n}{n^{2+2r-2s}}$,
then  almost surely the empirical distribution of eigenvalues
of the $s$-th Laplacian of  $H^r(n,p)$
 follows the Semicircle Law centered at $1$ and   with radius $(2+o(1))\sqrt{\frac{1-p}{{r-s\choose s}{n-s\choose r-s}p}}$.
 \end{theorem}

 \begin{remark}
   The proof of Theorem \ref{t2} actually implies the eigenvalues of
   $\L^{(s)}(H^r(n,p))-\L^{(s)}(K^r_n)$ follows the Semicircle Law
   centered at $0$ and with radius
   $(2+o(1))\sqrt{\frac{1-p}{{r-s\choose s}{n-s\choose r-s}p}}$. Thus
   we have
\begin{equation}
  \label{eq:difflb}
\max_{1\leq k \leq {n\choose s}-1}|\lambda_k^{(s)}(H^r(n,p))-\lambda_k^{(s)}(K^r_n)|\geq \left(\frac{2}{\sqrt{{r-s\choose s}}}+o(1)\right)\sqrt{\frac{1-p}{{n-s\choose r-s}p}}.
\end{equation}
This shows that the upper bound of
$|\lambda_k^{(s)}(H^r(n,p))-\lambda_k^{(s)}(K^r_n)|$  in Theorem
\ref{t1}  is best up to a constant multiplicative factor.
 \end{remark}

The rest of the paper is organized as follows. In section 2, we introduce
the notation and prove some basic lemmas. We will prove Theorem 1 in section 3
and Theorem 2 in section 4.

\section{Notation and Lemmas}

\subsection{Laplacian eigenvalues of hypergraphs}
Let $H=(V,E)$ be an $r$-uniform hypergraph. For any subset $S$
($|S|<r$), the degree of $S$, denoted by $d_S$, is the number of edges
passing through $S$.
 For each $1\leq s \leq
r/2$, we associate a weighted graph $G^{(s)}$ on the vertex set
${V\choose s}$ to $H$ as follows. Every pair of $s$-sets $S$ and $T$ is
associated with a weight $w(S,T)$, which is given by
$$w(S,T)=\left\{
  \begin{array}[c]{ll}
    d_{S\cup T} & \mbox{ if } S\cap T=\emptyset,\\
   0 & \mbox{ otherwise }.
  \end{array}
\right.
$$
The $s$-th Laplacian $\L^{(s)}$ of $H$ is defined
to be the normalized Laplacian of $G^{(s)}$. The degree of $S$
in $G^{(s)}$ is $\sum_{T}w(S,T)={r-s \choose s} d_S$.

We  assume that the $s$-sets in ${V\choose s}$ are ordered alphabetically.
Let $N:={n\choose s}$; all square matrices considered in the paper have the
dimension $N\times N$ and all vectors have dimension $N$.
Let $W:=(W(S,T))$ be the weight matrix,
$D$ be the diagonal matrix with diagonal entries $D(S,S)=d_S$,
$\d$ be the column vector with entries $d_S$ at position $S\in {V\choose S}$,
 $J$ be the square matrix of all $1$'s, and $\1$ be the column vector of all $1$'s.
Let $T:={r-s \choose s} D$; here $T$ is the diagonal matrix of degrees in $G^{(s)}$.
Then, we have
$$\L^{(s)}=I-T^{-1/2}WT^{-1/2}.$$

We list the eigenvalues of $\L^{(s)}$ as $$0=\lambda_0^{(s)} \leq
\lambda_1^{(s)},\ldots,\lambda_{\binom{n}{s}-1}^{(s)} \leq 2.$$

 We aim to
compute  the spectral radius
$\bar\lambda^{(s)}(H)=\max_{i\not =0}|1-\lambda_i^{(s)}|$.
Let $\vol^{(s)}(H):=\sum_{S\in {V\choose s}}d_s$ and
$\phi_0:=\frac{1}{\sqrt{\vol^{(s)}(H)}}D^{1/2}\1$. Note that $\phi_0$ is
the unit eigenvector corresponding to the trivial eigenvalue $0$ of $\L^{(s)}$.

We are ready to prove theorem \ref{t0}.

{\bf Proof of Theorem \ref{t0}:} We can express $\L^{(s)}(K^r_n)$ using
the following notation.
The  Kneser graph $K(n, s)$ is a graph over the vertex set ${[n]\choose s}$;
two $s$-sets $S$ and $T$ form an edge of $K(n, s)$ if and only if $S\cap T=0$.
Let $K$ be the adjacency matrix of $K(n,s)$; the eigenvalues of $K$ are
$(-1)^{i}\binom{n-s-i)}{s-i}$ with multiplicity ${n\choose i}-{n\choose i-1}$
for $0\leq i\leq s$ (see \cite{gr}). Note that $K(n, s)$ is a regular graph;
 so the Laplacian eigenvalues
can be determined from the eigenvalues of its adjacency matrix.
We observe that the  associated weighted graph
$G^{(s)}$ for the complete $r$-uniform hypergraph $K^r_n$ is essentially the Kneser graph
with each edge associated with a weight ${n-2s\choose r-2s}$.
Note that the multiplicative factor ${n-2s\choose r-2s}$ is canceled after normalization.
The $\L^{(s)}$ (for $K^r_n$) is
exactly
the Laplacian of Kneser graph.  Hence,
$$\L^{(s)}(K^r_n)=I-\frac{1}{{n-s\choose s}}K.$$
Thus, the eigenvalues of $s$-th Laplacian of $K^r_n$ are given by
$$
1- \frac{(-1)^{i}\binom{n-s-i)}{s-i}}{{n-s \choose s}}
\mbox{ with multiplicity } {n\choose i}-{n\choose i-1}
\mbox{ for }
0\leq i \leq s. $$
\hfill $\square$

\begin{remark}
 For $1\leq s\leq r/2$, we have
\begin{eqnarray}
  \label{eq:lambda1_knr}
\lambda_1^{(s)}(K^r_n) &=& 1- \frac{s(s-1)}{(n-s)(n-s-1)},\\
  \label{eq:lambdamax_knr}
\lambda_{max}^{(s)}(K^r_n)&=&1+\frac{s}{n-s},\\
  \label{eq:rho_knr}
\bar\lambda^{(s)}(K^r_n)&=&\frac{s}{n-s}.
\end{eqnarray}
  \end{remark}

\subsection{Random hypergraphs}
 Let $H^r(n,p)$ be a random $r$-uniform hypergraph over the vertex set $V=[n]$
 and each $r$-set has probability $p$ to be an edge independently.
We would like to bound the spectral radius of the $s$-th Laplacian
of $H^r(n,p)$ for $1\leq s\leq r/2$.

For any $F\in {V\choose r}$, let $X_F$ be the random indicator variable for
$F$ being an edge in $H^r(n,p)$; all $X_F$'s are independent to each other.
For any $S, T \in {V\choose s}$, we have
$$W(S,T)=\left \{
    \begin{array}{ll}
\sum_{\substack{{F\in {n\choose r}}\\ S\cup T\subset F}}X_F
       & \mbox{ if } S\cap T= \emptyset;\\
 0 & \mbox{ otherwise.}
    \end{array}
\right.
$$
Thus,
\begin{equation}
\label{eq:EW}
  \E(W(S,T))=\left \{
    \begin{array}{ll}
      \binom{n-2s}{r-2s}p & \mbox{ if } S\cap T= \emptyset;\\
 0 & \mbox{ otherwise.}
    \end{array}
\right.
\end{equation}
The degree $d_S=\sum_{S\subset F\in {V\choose r}}X_F$;
we have $\E(d_S)={n-s \choose r-s}p$. For simplicity, let
$d:={n-s \choose r-s}p$.

We use the following Lemma to compare the eigenvalues of two matrices.
\begin{lemma} \label{difference}
Given any two $(N\times N)$-Hermitian matrices $A$ and $B$, for $1\leq k\leq N$,
let $\mu_k (A)$ (or $\mu_k(B)$) be the $k$-th eigenvalues of $A$ (or $B$)
in the increasing order.
We have
$$|\mu_k(A)-\mu_k(B)|\leq \|A-B\|.$$
\end{lemma}
{\bf Proof:}  By the Min-Max Theorem (see \cite{minmax}), we have
\begin{eqnarray*}
\mu_k(A)&=&\min_{S_k}\max_{x\in S_k, \|x\|=1}x'Ax,\\
\mu_k(B)&=&\min_{S_k}\max_{x\in S_k, \|x\|=1}x'Bx.
\end{eqnarray*}
where the minimum is taken over all $k$-th dimensional subspace
$S_k\subset \mathbb{R}^N$.
We have
\begin{eqnarray*}
  \mu_k(A)&=&\min_{S_k}\max_{x\in S_k, \|x\|=1}x'Ax\\
&=&\min_{S_k}\max_{x\in S_k, \|x\|=1}(x'Bx+x'(A-B)x)\\
&\leq& \min_{S_k}\max_{x\in S_k, \|x\|=1}(x'Bx+\|A-B\|)\\
&=&\mu_k(B)+\|A-B\|.
\end{eqnarray*}
Similarly, we can show $\mu_k(A)\geq \mu_k(B)-\|A-B\|$.
The proof of the Lemma is finished. \hfill $\square$

Our idea is to bound the spectral norm of the difference of $\L^{(s)}(H^r(n,p))$
and $\L^{(s)}(K^r_n)$.
Let $M:=\L^{(s)}(K^r_n)-\L^{(s)}(H^r(n,p))=
T^{-1/2}WT^{-1/2}-\frac{1}{{n-s\choose s}}K$.
We write  $M=M_1+M_2+M_3+M_4$,
where
\begin{eqnarray*}
  M_1&=&\frac{1}{{r-s\choose s}}\left(D^{-1/2}(W-\E(W))D^{-1/2}-d^{-1}(W-\E(W))\right),
\\
  M_2 &=& \frac{1}{{r-s\choose s}d}(W-\E(W)),\\
  M_3 &=& \frac{1}{{r-s\choose s}}D^{-1/2}\E(W)D^{-1/2}-\frac{d}{{n\choose s}}D^{-1/2}JD^{-1/2}-\frac{1}{{n-s\choose s}}K+\frac{1}{{n\choose s}}J,\\
   M_4&=& \frac{1}{{n\choose s}}(dD^{-1/2}JD^{-1/2}-J).
\end{eqnarray*}
By the triangular inequality of matrix norms, we have
$$\|M\|\leq \|M_1\| +\|M_2\| + \|M_3\| +\|M_4\|.$$
Through this paper, the norm of any square matrix is the spectral norm.
We would like to bound $\|M_i\|$  for $i=1,2,3,4$.
We use the following Chernoff inequality.

\begin{theorem}{\cite{chernoff}}
\label{concen:chernoff}
 Let $X_1,\ldots,X_n$ be independent random variables with
$$\Pr(X_i=1)=p, \qquad \Pr(X_i=0)=1-p.$$
We consider the sum $X=\sum_{i=1}^n X_i$,
with expectation $\E(X)=np$. Then we have
\begin{eqnarray*}
\mbox{(Lower tail)~~~~~~~~~~~~~~~~~}
\qquad \qquad  \Pr(X \leq \E(X)-\lambda)&\leq& e^{-\lambda^2/2\E(X)},\\
\mbox{(Upper tail)~~~~~~~~~~~~~~~~~}
\qquad \qquad
\Pr(X \geq \E(X)+\lambda)&\leq& e^{-\frac{\lambda^2}{2(\E(X) + \lambda/3)}}.
\end{eqnarray*}
\end{theorem}

\begin{lemma}\label{degree}
Suppose $d \geq \log N$. With probability at least
$1-\frac{1}{N^2}$,
 we have $d_S\in (d-3\sqrt{d\log N}, d+3\sqrt{d\log N})$ for all $S\in {V\choose s}$.
\end{lemma}

{\bf Proof:} Note $d_s=\sum_{F: S\subset F}X_F$ and $\E(d_S)=d$.
Applying the lower tail of Chernoff's inequality with
$\lambda=3\sqrt{\E(X)\log N}$, we have
$$\Pr \left (X-\E(X)\leq - \lambda \right ) \leq e^{-\lambda^2/2\E(X)} \leq \frac{1}{N^{9/2}}.$$

Applying the upper tail of Chernoff's inequality with
$\lambda=3\sqrt{\E(X)\log N}$, we have
$$\Pr \left(X-\E(X) \geq  \lambda \right )
\leq e^{-\frac{\lambda^2}{2(\E(X) +
 \lambda/3)}}
\leq\frac{1}{N^{27/8}}.$$  
The probability that $d_S\not\in (d-3\sqrt{d\log N}, d+3\sqrt{d\log N})$
is at most $\frac{1}{N^3}$. Thus, with probability at least
$1-\frac{1}{N^2}$,
 we have $d_S\in (d-3\sqrt{d\log N}, d+3\sqrt{d\log N})$ for all $S\in {V\choose s}$.
\hfill $\square$

For convenience, let  $d_{min}:=d-3\sqrt{d\log N}$, $d_{max}:=d+3\sqrt{d\log N}$;
almost surely we have $d_{min}\leq d_S\leq d_{max}$ for all $S$.

\begin{lemma}
  \label{M3}
If $d\geq \log N$, then almost surely
$\|M_3\|=O\left( \frac{\sqrt{\log N}}{n\sqrt{d}}\right).$
\end{lemma}

{\bf Proof: }
Note  $\E(W)={n-2s\choose r-2s}pK$, where $K$ is
the adjacency matrix of the Kneser graph $K(n,s)$.
Let $M_0:=\frac{1}{{n-s\choose s}}K-\frac{1}{{n\choose s}}J$.
We can rewrite $M_3$ as
$$M_3= dD^{-1/2}M_0D^{-1/2}-M_0.$$
Note $\|M_0\|=\bar \lambda^{(s)}(K^r_n)=\frac{s}{n-s}$.
We have

\begin{eqnarray*}
  \|M_3\| &=&\|dD^{-1/2}M_0D^{-1/2}-M_0\|\\
&\leq & \|(dD^{-1/2}-d^{1/2}I)M_0D^{-1/2}\| +\|M_0( d^{1/2}D^{-1/2}-I)\|\\
&\leq & \|(d^{1/2}I-dD^{-1/2})\|\|M_0\|\|D^{-1/2}\| +\|M_0\|\|( d^{1/2}D^{-1/2}-I)\|\\
&\leq & \left |d^{1/2}-d d_{min}^{-1/2}\right|
\frac{s}{n-s} d_{min}^{-1/2} + \frac{s}{n-s}\left |d^{1/2}d_{min}^{-1/2} -1\right| \\
&=&O\left( \frac{\sqrt{\log N}}{n\sqrt{d}}\right).
\end{eqnarray*}
\hfill $\square$

\begin{lemma} \label{lemma5}
If $p(1-p) \gg \frac{\log n}{n^{r-s}}$, then  almost surely
$$\sum_{S\in {V\choose s}} (d_S-d)^2=(1+o(1)){n\choose s}d (1-p).$$
\end{lemma}

{\bf Proof:} For $S \in \binom{V}{s}$, let $X_S=(d_S-d)^2$. We have
$$\E(X_S)=\E((d_S-d)^2)=\Var(d_S)={n-s\choose r-s}p(1-p)=d(1-p).$$
We use the second moment method to prove that $\sum_SX_s$ concentrates
around its expectation ${n\choose s}d (1-p)$.
For any $S,T\in {V \choose s}$, the covariance can be calculated
as follows.

\begin{eqnarray*}
  \Cov(X_S, X_T)&=&
\E(X_SX_T)-\E(X_S)\E(X_T)\\
&=&\E((d_S-d)^2(d_T-d)^2)- d^2(1-p)^2.
\end{eqnarray*}

For $F \in \binom{V}{r}$, let $Y_F=X_F-\E(X_F)$. Then we have $d_S-d=\sum_{S\subset F}Y_F$.
$$\E((d_S-d)^2(d_T-d)^2)=\sum_{\substack{F_1,F_2\colon S\subset F_1\cap F_2\\
F_3,F_4\colon T\subset F_3\cap F_4}} \E(Y_{F_1}Y_{F_2}Y_{F_3}Y_{F_4}).$$
Since $\E(Y_{F_i})=0$, the non-zero terms occur only if
\begin{enumerate}
 \item $F_1=F_2=F_3=F_4$. In this case, we have
$$\E(Y_{F_1}Y_{F_2}Y_{F_3}Y_{F_4})=\E(Y_{F_1}^4)=(1-p)^4p+(-p)^4(1-p)=p(1-p)(1-3p+3p^2).$$
The number of choices is ${n-|S\cup T| \choose r-|S\cup T|}$.
\item  $F_1=F_2\not=F_3=F_4$. In this case, we have
$$\E(Y_{F_1}Y_{F_2}Y_{F_3}Y_{F_4})=\E(Y_{F_1}^2)\E(Y_{F_3}^2)=p^2(1-p)^2.$$
The number of choices is ${n-s\choose r-|S|}{n-s\choose r-|T|}- {n-|S\cup T|\choose r-|S\cup T|}$.
\item  $F_1=F_3\not=F_2=F_4$. In this case, we have
$$\E(Y_{F_1}Y_{F_2}Y_{F_3}Y_{F_4})=\E(Y_{F_1}^2)\E(Y_{F_2}^2)=p^2(1-p)^2.$$
The number of choices is ${n-|S\cup T|\choose r-|S\cup T|}^2-{n-|S\cup T|\choose r-|S\cup T|}$.
\item  $F_1=F_4\not=F_2=F_3$. This is the same as item 3.
\end{enumerate}
Thus, we have
\begin{eqnarray*}
\E(X_SX_T)&=&{n-|S\cup T|\choose r-|S\cup T|}p(1-p)(1-3p+3p^2)\\
&&+\left({n-s\choose r-s}^2 +2{n-|S\cup T|\choose r-|S\cup T|}^2
- 3{n-|S\cup T|\choose r-|S\cup T|}\right)p^2(1-p)^2\\
&=&{n-|S\cup T|\choose r-|S\cup T|}p(1-p)(1-6p+6p^2) +\left({n-s\choose r-s}^2 +2{n-|S\cup T|\choose r-|S\cup T|}^2\right)p^2(1-p)^2.
\end{eqnarray*}
This expression on the right depends only on the size of $S\cup T$.
Putting together, we get
\begin{eqnarray*}
  \Var\left(\sum_{S\in {V \choose s}}X_S\right)&=& \sum_{S,T\in  {V \choose s}}\Cov(X_S,X_T)\\
&=&\sum_{S,T\in  {V \choose s}}(\E(X_SX_T)-d^2(1-p)^2)\\
&=&\sum_{S,T\in  {V \choose s}}\left(\E(X_SX_T)-{n-s \choose r-s}^2p^2(1-p)^2\right)\\
&=& \sum_{i=s}^{2s}\sum_{|S\cup T|=i}
\left({n-i\choose r-i}p(1-p)(1-6p+6p^2)
+2{n-i\choose r-i}^2p^2(1-p)^2\right)\\
&\leq&  \sum_{i=s}^{2s}\sum_{|S\cup T|=i}
{n-i\choose r-i}p(1-p)\left(1-6p+6p^2+ 2{n-s\choose r-s} p(1-p)\right)\\
&\leq&  \sum_{i=s}^{2s}\sum_{|S\cup T|=i}
{n-i\choose r-i}3dp(1-p)^2\\
&=& {n\choose r}3dp(1-p)^2 \sum_{i=s}^{2s}\frac{r!}{(i-s)!^2(2s-i)!(r-i)!}\\
&<& 3\cdot 4^r  {n\choose r}dp(1-p)^2 \\
&=&O\left({n\choose s } d^2(1-p)^2\right ).
\end{eqnarray*}

Let $X=\sum_S X_S$. We have
 $\E[X]=\binom{n}{s}d(1-p)$ and $\Var(X) = O\left({n\choose s } d^2(1-p)^2\right )$.
Applying Chebyshev's inequality to $X=\sum_{S\in {V\choose s}}$,
we  have
$$
\Pr\left(|X-\E(X)| \geq \log n \sqrt{\Var(X)}\right) \leq \frac{1}{\log^2 n}.
$$
Thus, almost surely $X=\E(X)+O(\log n \sqrt{\Var(X)})=(1+o(1)) \binom{n}{s}d(1-p)$.  \hfill $\square$

\begin{lemma}\label{M4}
  If $p(1-p)\gg \frac{\log n}{n^{r-s}}$,
 then  almost surely
$\|M_4\|\leq(1+o(1))\sqrt{\frac{1-p}{d}}.$
\end{lemma}
{\bf Proof:}
We can rewrite $M_4$ as
\begin{eqnarray*}
M_4&=& \frac{1}{{n\choose s}}(
dD^{-1/2}JD^{-1/2}- J) \\
&=& \frac{1}{{n\choose s}} \left(
\left(d^{1/2}D^{-1/2} -I\right) JD^{-1/2}d^{1/2}+ J \left(d^{1/2}D^{-1/2} -I\right)
\right)\\
&=&  \frac{1}{{n\choose s}}\left(\alpha \1'D^{-1/2}d^{1/2}+   \1\alpha'\right).
\end{eqnarray*}
Here $\alpha:=d^{1/2}D^{-1/2}\1 -\1$.
Note that the spectral norm of a vector is the same as the $L_2$-norm.
We have
\begin{eqnarray*}
\|\alpha\|&=&  \|d^{1/2}D^{-1/2}\1 -\1\|\\
&=&\sqrt{\sum_{S\in {V\choose s}}\left(\frac{\sqrt{d}}{\sqrt{d_S}}-1\right)^2}\\
&=&\sqrt{\sum_{S\in {V\choose s}} \frac{(d_S-d)^2}{d_S(\sqrt{d}+\sqrt{d_S})^2}}
\\
&\leq & \frac{\sqrt{\sum_{S\in {V\choose s}}(d_S-d)^2}}{\sqrt{d_{min}}(\sqrt{d}+\sqrt{d_{min}})}\\
&=&(\frac{1}{2}+o(1))\sqrt{\frac{(1-p){n\choose s}}{d}}.
\end{eqnarray*}
In the last step, we applied Lemma \ref{lemma5}.
Therefore,
we have
\begin{eqnarray*}
 \|M_4\| &=& \left \|  \frac{1}{{n\choose s}}\left(\alpha \1'D^{-1/2}d^{1/2}+   \1\alpha'\right)
\right\|\\
&=&
 \frac{1}{{n\choose s}}\left(\left \|\alpha \1'D^{-1/2}d^{1/2}\right\|
+ \left\| \1\alpha' \right\|\right)
\\
&\leq&  \frac{1}{{n\choose s}} \|\alpha\| \left(\|\1'D^{-1/2}d^{1/2}\|
+\|\1\|\right)\\
&=&  \frac{1}{{n\choose s}} \|\alpha\| \left(
\sqrt{\sum_{S\in{n\choose s}}
\frac{d}{d_S}}+ \sqrt{{n\choose s}}\right)\\
&\leq&  \frac{1}{{n\choose s}}
\left(\frac{1}{2}+o(1)\right)\sqrt{\frac{(1-p){n\choose s}}{d}}
(2+o(1))\left(\sqrt{{n\choose s}} \right)\\
&=& (1+o(1))\sqrt{\frac{1-p}{d}}.
\end{eqnarray*}

\section{Proof of Theorem \ref{t1}}
To estimate the spectral norm of $M_1$ and $M_2$, we need consider
the matrix $C:=W-\E(W)$. We estimate the expectation of  the trace
of $C^t$ as follows.

\begin{lemma}\label{trace}
For any $k$ satisfying $k\ll \sqrt[4]{n^{r-s} p(1-p)}$, we have
\begin{eqnarray}
  \label{eq:enentrace2}
  \E\left(\tr(C^{2k})\right)&\leq &   (1+o(1))
\frac{n^{s+k(r-s)}{r-s \choose s}^k }{(k+1) (s!)^{k+1}((r-2s)!)^k} \binom{2k}{k} p^k(1-p)^k,\\
  \label{eq:oldtrace}
  \E\left(\tr(C^{2k+1})\right)&=&  O\left(\frac{2k^2n^{s+k(r-s)}{r-s\choose s}^{k+1}}{(k+1)
(s!)^{k+1} ((r-2s)!)^k} \binom{2k}{k} p^k(1-p)^k\right).
\end{eqnarray}

If  further $k =o\left(\log (n^{r-s}p(1-p))\right)$, then we have
\begin{eqnarray}
  \label{eq:eventrace}
\E\left(\tr(C^{2k})\right)&=&   (1+o(1))
\frac{n^{s+k(r-s)}}{(k+1) (s!)^{k+1} ((r-2s)!)^k} \binom{2k}{k} p^k(1-p)^k.
\end{eqnarray}
\end{lemma}
The proof of this technical Lemma is quite long. We will delay its proof
until the end of this section.

\begin{lemma}\label{C}
If $p(1-p)\gg \frac{\log^4 n}{n^{r-s}}$, then  we have $\|C\|\leq
\left(2{r-s\choose s}+o(1)\right)\sqrt{d(1-p)}$ almost surely.
\end{lemma}
{\bf Proof:} By Lemma \ref{trace}, we have $\E(\tr(C^{2k})) \leq
(1+o(1)) \frac{n^{s+k(r-s)}{r-s\choose s}^k}{(k+1) (s!)^{k+1}
((r-2s)!)^k}  \binom{2k}{k} p^k(1-p)^k$. As $\E(\|C\|^{2k}) \leq
\E(\tr(C^{2k}))$, we have
$$
\E(\|C\|^{2k}) \leq
(1+o(1))\frac{n^{s+k(r-s)}\binom{r-s}{s}^k}{(k+1) (s!)^{k+1}
((r-2s)!)^k} \binom{2k}{k} p^k (1-p)^k.
$$
Let $U:=\frac{n^{s+k(r-s)}{r-s\choose s}^k}{(k+1) (s!)^{k+1}
((r-2s)!)^k} \binom{2k}{k} p^k(1-p)^k$.  By Markov's inequality,
\begin{eqnarray*}
\Pr  \left(\|C\| \geq (1+\epsilon)\sqrt[2k]{U} \right)
&=& \Pr\left(\|C\|^{2k} \geq (1+\epsilon)^{2k}U\right) \\
& \leq & \frac{\E(\|C\|^{2k})}{(1+\epsilon)^{2k}U} \\
&\leq&
\frac{(1+o(1)) U}{(1+\epsilon)^{2k}U} \\
&=& \frac{1+o(1)}{(1+\epsilon)^{2k}}.
\end{eqnarray*}
Let $g(n)$ be a slowly growing function such that $g(n)\to \infty$ as $n$ approaches
the infinity and $g(n)\ll \frac{(n^{r-s}p(1-p))^{1/4}}{s \log n}$. This is possible
because  $n^{r-s}p(1-p)   \gg \log^{4} n$.
Choose  $k= sg(n)\log n$ and $\epsilon=1/g(n)$.
We have $k\ll (n^{r-s}p(1-p))^{1/4}$ and $\epsilon\to 0$.
 Then
 we have $(1+o(1))/(1+\epsilon)^{2k}=O(n^{-s})$, which implies that almost surely
  \begin{eqnarray*}
|\|C\| &\leq& (1+o(1)) \sqrt[2k]{U}\\
&=& (1+o(1)) \left( \frac{n^{s+k(r-s)}\binom{r-s}{s}^k}{(k+1) (s!)^{k+1} ((r-2s)!)^k} \binom{2k}{k}p^k(1-p)^k\right)^{\frac{1}{2k}}\\
&<&  n^{\frac{s}{2k}}2\sqrt{\frac{n^{r-s}\binom{r-s}{s}p(1-p)}{s!(r-2s)!}}\\
&=&\left(2\binom{r-s}{s}+o(1)\right) \sqrt{d(1-p)}.
 \end{eqnarray*}
 \hfill
$\square$

Recall  $M_2=\frac{1}{{r-s\choose s}d}C$. We have
\begin{lemma}\label{M2}
If $p(1-p)\gg \frac{\log^4 n}{n^{r-s}}$, then we have $\|M_2\|\leq
(2+o(1))\sqrt{\frac{1-p}{d}}$ almost surely.
\end{lemma}

\begin{lemma} \label{M1}
  If $p(1-p)\gg \frac{\log^4 n}{n^{r-s}}$, then   we have
$\|M_1\|=O\left(\frac{\sqrt{(1-p)\log N}}{d}\right)$ almost surely.
\end{lemma}
{\bf Proof:}
We have
\begin{eqnarray*}
 M_1&=&\frac{1}{{r-s\choose s}}\left(D^{-1/2}CD^{-1/2}-d^{-1}C\right)\\
 &=&\frac{1}{{r-s\choose s}} \left((D^{-1/2}-d^{-1/2}I) C D^{-1/2}
+d^{-1/2}C(D^{-1/2}-d^{-1/2}I)\right).
\end{eqnarray*}

Note $\|D^{-1/2}-d^{-1/2}I\|\leq
|d_{min}^{-1/2}-d^{-1/2}|=O(\frac{\sqrt{\log N}}{d})$,
$\|D^{-1/2}\|\leq d_{min}^{-1/2}=(1+o(1))d^{-1/2}$, and
$\|C\|=\left(2\binom{r-s}{s}+o(1)\right) \sqrt{d(1-p)}$. We have
\begin{eqnarray*}
 \|M_1\| &=&\frac{1}{{r-s\choose s}} \left\| (D^{-1/2}-d^{-1/2}I) C D^{-1/2}
+d^{-1/2}C(D^{-1/2}-d^{-1/2}I)\right\|\\
&=&O\left(\frac{\sqrt{(1-p)\log N}}{d}\right).
\end{eqnarray*}
 \hfill
$\square$

{\bf Proof of Theorem \ref{t1}:}
Combining  Lemmas \ref{M3}, \ref{M4}, \ref{M2}, and \ref{M1},
we have

\begin{eqnarray*}
 \|M\|&=& \|M_1+M_2+M_3+M_4\|\\
&\leq&  \|M_1\|+\|M_2\|+\|M_3\|+\|M_4\|\\
&\leq & O\left(\frac{\sqrt{(1-p)\log N}}{d}\right)
+ \frac{(2+o(1))\sqrt{1-p}}{\sqrt{d}}+ O\left(\frac{\sqrt{\log N}}
{n\sqrt{d}}\right)
+ (1+o(1)) \sqrt{\frac{1-p}{d}}\\
&=& \left(3+o(1)\right)\sqrt{\frac{1-p}{d}}.
\end{eqnarray*}
In the last step, we use the fact $\frac{\sqrt{\log N}}{n\sqrt{d}}
=o\left(\sqrt{\frac{1-p}{d}}\right)$ since $1-p\gg \frac{\log n}{n^2}$.

By Lemma \ref{difference}, for $1\leq k \leq {n\choose s}-1$, we have
$$|\lambda^{(s)}_k(H^r(n,p))-\lambda^{(s)}_k(K^r_n)|\leq
\|M\|\leq  \left(3+o(1)\right)\sqrt{\frac{1-p}{d}}.$$
\hfill
$\square$

Reall that $X_F$ is the random indicator variable for
$F$ being an edge in $H^r(n,p)$.
For any fixed positive integer $t$,  the terms in $\tr(C^{t})$ are of the form
$$
c_{S_1S_2}c_{S_2S_3}\ldots c_{S_{t}S_{S_1}}.
$$
Here $c_{ST}=W(S,T)-\E(W(S,T))=\sum_{\substack{F\in {V\choose r}\\
S\cup T\subset F}}(X_F -\E(X_F))$ if $S\cap T=\emptyset$; $c_{ST}=0$ otherwise.

Note   $c_{S_iS_{j}}=0$ if $S_i \cap S_j \not = \emptyset$. Thus
we need only to consider the sequence $S_1S_2\ldots S_{t}S_1$ such
that $S_{i} \cap S_{i+1}=\emptyset$ for each $1 \leq i \leq t$, here
$t+1=1.$

For  $F\in {V\choose r}$ and $S,T\in {V\choose s}$,
 we define a  random variable $c_{ST}^F$ as follows.
$$
 c_{ST}^F = \left\{
  \begin{array}{ll}
      X_F-\E(X_F) &  \mbox{ if }  S\cap T=\emptyset \mbox{ and }
S\cup  T \subseteq F;\\
    0 &  \mbox{ otherwise.}
  \end{array}
   \right.
$$
The sequence  $w:=S_1F_1S_2F_2S_3 \ldots
S_{t} F_{t}S_{1}$ is called a closed $s$-walk of length $t$ if
\begin{enumerate}
\item $S_1,\ldots, S_t\in {V\choose s}$,
\item $F_1,\ldots, F_t \in {V\choose r}$,
\item $S_i\cap S_{i+1}=\emptyset$, for $i=1,2,\ldots, t$,
\item $S_i\cup S_{i+1}\subset F_i$, for $i=1,2,\ldots, t$.
\end{enumerate}
Here we use the convention $S_{t+1}=S_1$. Those $r$-sets $F_i$'s are referred
as edges while those $s$-sets $S_i$'s are referred as stops.
For $1\leq i \leq t$, we say $w$ walks from $S_i$ to $S_{i+1}$ at step $i$
via the edge $F_i$.

Using the notation above, we rewrite the trace  as
$$
\tr(C^{t})= \sum_{\mbox{closed $s$-walks}} c_{S_1S_2}^{F_1} c_{S_2S_3}^{F_2} \ldots
c_{S_{t}S_1}^{F_{t}},
$$
where the summation is over all possible closed $s$-walks
of length $t$.

Taking the expectation on both sides, we get
$$
\E(\tr(C^{t}))= \sum_{\mbox{closed $s$-walks}}  \E(c_{S_1S_2}^{F_1} c_{S_2S_3}^{F_2} \ldots
c_{S_{t}S_1}^{F_{t}}).
$$

The terms in the product above can be regrouped according to the
values of $F_i$'s; those terms with distinct $F$'s are independent
to each other. Since $\E(c^F_{S,T})=0$, the contribution of a closed
walk is $0$ if some $F$ appears only once.  Thus we need only to
consider the set of closed walks where each edge appears at least
twice or do not occur; we call these closed walks as {\em good}
closed walks.  A good closed walk can contain at most $\lfloor
\frac{t}{2}\rfloor$ distinct edges.

Let $\G_i$ be the set of good closed walks of length $t$ with
$i$ distinct edges. For $1\leq i \leq \lfloor \frac{t}{2}\rfloor$,
let ${\cal G}_i^j$ be the set of good closed walks with exactly $i$
distinct edges and $j$ distinct vertices; we have
$\G_i:=\cup_{j}\G_i^j$.

We consider  a  good closed walk  in ${\cal G}_i$. When a new edge
comes in the walk, it can bring in  at most $(r-s)$ new vertices.
Thus such a good closed walk covers at most $m_i:=s+i(r-s)$
vertices. Any walk contains at least one edge. Hence, the number of
vertices in a walk from $\G_i$ is in the interval $[r, m_i]$.

We have
\begin{equation}\label{eq:trace}
\E(\tr(C^{t}))=
\sum_{i=1}^{\lfloor \frac{t}{2}\rfloor}
\sum_{ S_1F_1S_2\ldots S_{t}S_1\in {\cal G}_i} \E(c_{S_1S_2}^{F_1}
 c_{S_2S_3}^{F_2} \ldots c_{S_{t}S_1}^{F_{t}}).
\end{equation}

Assume that an edge $F$ occurs $q$ times in a good closed walk
 and $T:=\{i: 1 \leq i \leq t \
\textrm{and} \ F_i=F\}$.  We have  $\Pr \left( \Pi_{i \in T}
c_{S_{i}S_{i+1}}^F=(1-p)^q\right)=p$ and $\Pr\left(\Pi_{i \in T}
c_{S_{i}S_{i+1}}^F=(-p)^q\right)= 1-p$. Thus, for each positive
integer $l \geq 2$, we have
$$
\E\left(\Pi_{i \in T} c_{S_{i}S_{i+1}}^F\right) =(1-p)^q
p+(-p)^q(1-p) \leq p(1-p).
$$
The equality holds for $q=2$.

Pick  a good closed walk $w:=S_1F_1S_2F_2S_3 \ldots
S_{t} F_{t}S_{1}$ in ${\cal G}_i$. Let $F^1,\ldots, F^i$ be the list of
distinct edges in the order as they appear in $w$. 

 For each $1 \leq l \leq i$, let $T_l:=\{1
\leq j \leq t:  F_j=F^l \}$; then $\sum_{l=1}^i |T_l|=t.$ We have
$$
\E(c_{S_1S_2}^{F_1} c_{S_2S_3}^{F_2} \ldots
c_{S_{t}S_1}^{F_{t}})=\Pi_{l=1}^i \Pi_{j \in T_l} \E(c_{S_j
S_{j+1}}^{F^l}) \leq \Pi_{l=1}^i p(1-p)=p^i (1-p)^i.
$$

This implies
\begin{equation}
  \label{eq:ai}
 \sum_{ S_1F_1S_2\ldots S_{t}S_1\in {\cal G}_i} \E(c_{S_1S_2}^{F_1}
 c_{S_2S_3}^{F_2} \ldots c_{S_{t}S_1}^{F_{t}}) \leq \left|{\cal G}_i\right| p^i(1-p)^i
\end{equation}
  for all $1\leq i\leq {\lfloor \frac{t}{2}\rfloor}$.
In particular, the equality holds when $t=2i$. Combining equation
(\ref{eq:trace}) and inequality (\ref{eq:ai}), we get
\begin{equation}
  \label{eq:trace1}
  \E(\tr(C^{t}))\leq
\sum_{i=1}^{\lfloor \frac{t}{2}\rfloor}|\G_i| p^i(1-p)^i.
\end{equation}

Now we estimate the value of $|\mathcal G_i^j|$, the number of good
closed walks of length $t$ on $i$ edges and $j$ vertices.  
Let $w$ be a good closed walk in $\G_i^j$.
For $2\leq k \leq
i$, let $\cdots SF^kS'\cdots$ be a piece of sequence in $w$ where
the edge $F^k$ occurs first time; $S$ is called the in-stop of $F^k$ and $S'$ is
called the out-stop of $F^k$.

The following lemma will state the hypergraph structure of these $i$
edges; it is independent of the walk $w$. We will use the following
notation.  Let $\S=\cup_{l=1}^i{F^l\choose s}$.  For any $s$-set
$S\in \S$, the degree of $S$, denoted by $d_S$, is the number of
edges in $\{F^1,F^2,\ldots, F^i\}$ containing $S$.

Define
$$d'_S=\left\{
  \begin{array}{ll}
    d_S-1 & \mbox { if there exists a unique $k$ such that } S=F_k\cap (\cup_{l=1}^{k-1}F_l),\\
d_S & \mbox{ otherwise}.
  \end{array}
\right.
$$


\begin{lemma}\label{l:iedges}
Assume that $F^1,\ldots, F^i$ is the list of distinct edges in the
order as they appear in $w\in \G_i^j$. Then we have
$$\sum_{S\in \S} (d'_S -1)\leq
\left(1 +\frac{2}{s}{r\choose s-1}\right ) (m_i-j).$$
\end{lemma}

{\bf Proof:}
For $2\leq k\leq i$,  let $x_k=|F^{k}\setminus (\cup_{l=1}^{k-1}F^l)|$;
we have $$0\leq x_k\leq r-s.$$
Thus,
$$j= r+x_2+x_3 +\cdots +x_{i}\leq r+(i-1)(r-s)=m_i.$$
Since a new edge $F_k$ can  contribute at most ${r-x_k\choose s}$ to
$\sum_{S\in \S} (d_S -1)$, we have
$$\sum_{S\in \S} (d_S -1)\leq \sum_{k=2}^i {r-x_k\choose s}.$$
Let $K:=\{k\colon x_k=r-s, 2\leq k\leq i\}$ and $\overline {K}=
\{2,\ldots, i\}\setminus K$. The edges in the set $\{F_k\colon k\in
K\}$ are called {\it  forward} edges while the edges in the set
$\{F_k\colon k\in \overline K\}$ are called {\it backward} edges.
Note each backward edge contribute at least one  to $m_i-j$; thus
$$m_i-j\geq |\overline {K}|.$$

Note for each $k\in K$,  ${r-x_k\choose s}=1$.
We have
\begin{eqnarray*}
\sum_{S\in \S} (d_S -1) &\leq &\sum_{k=2}^i {r-x_k\choose s}\\
&=& |K|+ \sum_{k\in \overline {K}}{r-x_k\choose s}\\
&=&|K|+ \sum_{k\in \overline {K}}\frac{r-x_k-s+1}{s} {r-x_k\choose s-1}\\
&\leq& |K|+ \sum_{k\in \overline {K}}\frac{2}{s}{r\choose s-1}(r-x_k-s)\\
&=& |K|+ \frac{2}{s}{r\choose s-1}(m_i-j).
\end{eqnarray*}

For any $k\in K$, let $S(F_k):=F_k\cap (\cup_{l=1}^{k-1}F_l)$ be the
starting stop of $F_k$  when $F_k$ first occurs in $w$. List the
elements in $K$ as $k_1,k_2,\ldots, k_{|K|}$ in an increasing order.
Consider the sequence of stops $S(F_{k_1}), S(F_{k_2}),\ldots,
S(F_{k_{|K|}})$ (not necessarily distinct). Let  $z$ be the number
of distinct stops in the sequence. If $S(F_{k_l})$ does not appear
the first time in the sequence above, then we consider the partial
walk $S_{k_{l-1}}F_{k_{l-1}} \ldots S_{k_l}F_{k_l}$. Since $F_{k_{l-1}}$ is a
forward edge,  there exists at
least one backward edge $F^{l'}$ for some $l'\in (k_{l-1}, k_{l})$.
 Thus,
$$|K|\leq z+|\overline {K}|\leq z +m_i-j.$$
Hence,
\begin{eqnarray*}
\sum_{S\in\S}(d'_S-1) &=& \sum_{S\in\S}(d_S-1) -z \\
&\leq&   |K|+ \frac{2}{s}{r\choose s-1}(m_i-j) -z\\
&\leq& \frac{2}{s}{r\choose s-1}(m_i-j) + m_i-j\\
&=&\left(1+\frac{2}{s}{r\choose s-1}\right) (m_i-j).
\end{eqnarray*}
The proof of this Lemma is finished. \hfill $\square$

\begin{lemma} \label{1:gij}
For $1\leq i \leq \lfloor \frac{t}{2} \rfloor$ and
$r\leq j \leq m_i$, we have
$$|\G_i^j| \leq {t-2\choose t-2i} i^{t-2i}\frac{1}{i+1}{2i\choose i} {r-s\choose s}^{t-i}
\frac{n^{m_i}}{(s!)^{i+1}((r-2s)!)^{i}} \left( \frac{C_1i^{C_2}}{n}\right)^{m_i-j}.
$$
Here $C_1$ and $C_2$ depend only on $r$ and $s$, independent of $i$,
$j$, and $n$.
\end{lemma}

\begin{corollary}
  For $1\leq i \leq \lfloor \frac{t}{2} \rfloor$ and $n\gg i^{C_2}$, we have
  \begin{equation}
    \label{eq:gi}
|\G_i| \leq (1+o(1)){t-2\choose t-2i} i^{t-2i}\frac{1}{i+1}{2i\choose i} {r-s\choose s}^{t-i}
\frac{n^{m_i}}{(s!)^{i+1}((r-2s)!)^{i}}.
  \end{equation}
\end{corollary}

{\bf Proof:}
We can associate a walk $w\in \G_i^j$ with a code of length $t$ consisting
of three symbols: `$($',  `$)$', and `$\ast$'.
 We scan the edges of the walk $w$ from
left to right; if an edge appears first time, then we assign the
code `$($'; if an edge appears second time,  then we assign the code
`$)$'; otherwise, we assign the code `$\ast$'.

For example, consider the following good walk
 with $i=3$, $j=8$, and $t=8$:
$$w=S_1 F_1 S_2 F_2 S_3 F_3 S_4  F_1 S_5  F_1 S_4 F_3 S_3 F_2 S_2 F_1 S_1$$
Here edges are: $F_1=(1,2,3,4,5)$, $F_2=(4,5,6,7,8)$,
$F_3=(7,8,5,3,4)$. Stops are: $S_1=(1,2)$, $S_2=(4,5)$, $S_3=(7,8)$,
$S_4=(3,,4)$, $S_5=(2,5)$. The code for this walk is $((()*))*$.

Since $w$ has $i$ distinct edges,  there are $i$ `$($'s,  $i$
`$)$'s, and  $(t-2i)$ `$\ast$'s. Note that the number of  `$($'  is
always greater than or equal to the number of  `$)$' at any point
when the sequence is read from left to right; each `$($' has a
matched `$)$' in the sequence.
 The symbol '$\ast$' starts at position three and up.
There are ${t-2\choose t-2i}$ ways to choose the '$\ast$'-positions
and $\frac{1}{i+1}{2i\choose i}$ ways to choose $i$ matched
parentheses (the Catalan number). The number of such codes is
$${t-2 \choose t-2i} \frac{1}{i+1}{2i\choose i}.$$

To construct a walk from a given code, we scan the symbols from left
to right. The first symbol is always `$($'. There are ${n\choose s}$
ways to choose the first stop $S_1$ and ${n-s \choose r-s}$ ways to
choose the rest of vertices in the first edge $F^1$. Suppose that we
already build a partial walk  and need to decide the next stop and
the next edge. There are at most ${r-s\choose s}$ ways to choose the
next stop $S$. The choices of selecting the next edge depends on the
next available symbol in the code sequence. Let $b_($, $b_)$, and
$b_\ast$ be the product of the number of ways to choose the next
edge at the `(', `)', and `$\ast$' positions respectively. We have
\begin{equation}
  \label{eq:gij0}
|\G_i^j|\leq {t-2 \choose t-2i} \frac{1}{i+1}{2i\choose i}{n\choose s} {r-s\choose s}^{t} b_( \cdot b_\ast \cdot b_).
\end{equation}



First we estimate $b_($, the number of ways to choose  new edges
$F^1,\ldots, F^i$ given the first stop $S_1$. Besides the $s$
vertices selected at the first stop, there are  ${n-s\choose j-s}$
ways to choose remaining $j-s$ vertices. Recall  that $F^1,\ldots,
F^i$ is the list of distinct edges in the order as they appear in
$w$. For $2\leq l\leq i$, let  $\tilde F^l:= F^l\setminus
(\cup_{l'=1}^{l-1}F^{l'})$,
 $x_l:=|\tilde F^l|$, and $y_l:=r-s -x_l$. We also define $\tilde F^1:=F^1\setminus S_1$;
$x_1:=|\tilde F^1|=r-s$, and $y_1=0$. Note that
$\cup_{l=1}^{i}\tilde F^l$ forms a partition of the remaining
$(j-s)$ selected vertices. The number of ways to choose such a
partition is
$$\frac{(j-s)!}{x_1! x_2!\cdots x_{i}!}.$$ To choose $F_l$, we need
select $x_l$ new vertices and $y_l$ old vertices; each old vertex
has at most $j$ choices. We have
$$b_(\leq \sum_{x_2,\ldots, x_i}{n-s\choose j-s} \frac{(j-s)!}{(r-s)! x_2!\cdots x_{i}!}
j^{\sum_{l=2}^i y_l}.$$
Observe  that $\sum_{l=2}^{i}y_l=m_i-j$ and
$$\frac{(j-s)!}{(r-s)! x_2!\cdots x_{i}!}\leq {j -s\choose r-s}\frac{(m_i-r)!}{((r-s)!)^{i-1}}.$$
The number of ways to choose $x_2,\ldots, x_{i}$ is the same as
the number of ways to choose $y_2,\ldots, y_{i}$, which
is   ${m_i-j+i-2\choose m_i-j}\leq (m_i-j+i-2)^{m_i-j}$. Therefore,
\begin{eqnarray}
\nonumber
b_( &\leq&  \sum_{x_2,\ldots, x_i}{n-s\choose j-s} \frac{(j-s)!}{(r-s)! x_2!\cdots x_{i}!}
j^{\sum_{l=2}^i y_l}\\
\nonumber
&\leq&  {n-s\choose j-s}{j-s\choose r-s} \frac{(m_i-r)!}{((r-s)!)^{i-1}} (m_i-j+i-2)^{m_i-j}
 j^{m_i-j}\\
&\leq&  {n-s\choose j-s}{j-s\choose r-s} \frac{(m_i-r)!}{((r-s)!)^{i-1}}
\left(\frac{m_i+i-2}{2}\right)^{2(m_i-j)}.
\label{eq:blp}
\end{eqnarray}

There is at most $i$ choices of edges at each `$\ast$' position. Thus
\begin{equation}
  \label{eq:bast}
 b_*\leq  i^{t-2i}.
\end{equation}

It remains to bound $b_)$. We first present an easy bound for $b_)$.
Edge $F$ can be chosen at most one  `$)$'-position. For any possible
stop $S\in \S$ , $S$ can appear at the $)$-positions at most $d_S$
times; each occurrence of $S$ involves different edges since we are
considering the second occurrence of edges. Thus,
$$b_)\leq \prod_{S\in \S} d_S! \leq  \prod_{S\in \S} d_S^{d_S-1}
\leq i^{\sum_{S\in \S} (d_S-1)}.$$ 

We need a better upper bound for $b_)$. Consider  a stop $S$ which  is first chosen
at  a `$)$'-position. Let $F$ be the edge on the walk right before
the `$)$'-position; i.e., the walk $w$ enters $S$ through $F$. 
If this $F$ occurred before, then the choices of edges at
`$)$'-positions starting with $S$ is at most
$$(d_S-1)!\leq i^{d_S-2}.$$
If this $F$ occurs first time and $F$ is an forward edge, 
then there is only one choice for the next edge leaving $S$;
namely $F$ itself. In this case, the choices of edges at
`$)$'-positions starting with $S$ is at most
$$(d_S-1)!\leq i^{d_S-2}.$$
In the remaining case, $F$ must be a backward edge.
The number of backward edges is at most $m_i-j$.
Since $F$ contains at most ${r\choose s}$ stops, the number of  such
$S$ is at most ${r\choose
  s}(m_i-j)$. A additional factor  $i^{{r\choose s}(m_i-j)}$ is enough.
 We have
 \begin{eqnarray}
\nonumber
b_)&\leq&  i^{{r\choose s}(m_i-j)} \prod_{S: d_S\geq 2} i^{d_S-2}\\
\nonumber
&\leq&  i^{{r\choose s}(m_i-j)} \prod_{S} i^{d'_S-1}\\
\nonumber
&\leq&    i^{{r\choose s}(m_i-j)}  i^{(1 +\frac{2}{s}{r\choose s-1})(m_i-j)}\\
&=& i^{\left ({r\choose s}+ 1 +\frac{2}{s}{r\choose s-1}\right)(m_i-j)}.
   \label{eq:brp}
 \end{eqnarray}

Combining equations (\ref{eq:gij0}), (\ref{eq:blp}),  (\ref{eq:bast}),  and (\ref{eq:brp}),
we get
\begin{eqnarray*}
|\G_i^j| &\leq&   {t-2 \choose t-2i} \frac{1}{i+1}{2i\choose i}{n\choose s} {r-s\choose s}^t {n-s\choose j-s}{j-s\choose r-s}
\\
&& \frac{(m_i-r)!}{((r-s)!)^{i-1}}
\left(\frac{m_i+i-2}{2}\right)^{2(m_i-j)}i^{t-2i}i^{\left({r\choose s}+ 1 +\frac{2}{s}{r\choose s-1}\right)(m_i-j)}\\
&\leq& {t-2\choose t-2i} i^{t-2i}\frac{1}{i+1}{2i\choose i} {r-s\choose s}^{t-i}
\frac{n^{j}}{(s!)^{i+1}((r-2s)!)^{i}} \\
&&\frac{(m_i-r)!}{(j-r)!}\left(\frac{m_i+i-2}{2}\right)^{2(m_i-j)}i^{\left({r\choose s}+ 1 +\frac{2}{s}{r\choose s-1}\right)(m_i-j)}\\
&\leq& {t-2\choose t-2i} i^{t-2i}\frac{1}{i+1}{2i\choose i} {r-s\choose s}^{t-i}
\frac{n^{m_i}}{(s!)^{i+1}((r-2s)!)^{i}} \left( \frac{C_1i^{C_2}}{n}\right)^{m_i-j}\\
\end{eqnarray*}
Here we set $C_1=4(r-s)^3$ and $C_2={r\choose s}+ 4 +\frac{2}{s}{r\choose s-1}$.
\hfill $\square$

\begin{lemma}
  If $t=2k$ is even, then we have
\begin{equation}
  \label{eq:gk}
|\G_k^{m_k}| ={n\choose m_k}\frac{m_k!}{(k+1)(s!)^{k+1}((r-2s)!)^{k}}\binom{2k}{k}.
\end{equation}
\end{lemma}
\noindent {\bf Proof:} We will construct a bijection from
$\G_k^{m_k}$ to a triple $(U, \P, C)$, where $U$ is a set of $m_k$
vertices, $\P$ is a partition of $U$ into $(k+1)$ $s$-sets and $k$
$(r-2s)$-sets, and $C$ is a code consisting of  $k$ pairs valid
parentheses.

For any good walk $w\in \G_k^{m_k}$,   let $U$ be the set of
vertices covered by $w$. Note each edge appears exactly twice. We
define a graph $T$, whose vertices are the  stops in $w$. Two stops
are connected if they belong to one edge. Observe that $T$ is
acyclic and connected; $T$ must be a tree. Since $T$ has exactly $k$
edges, $T$ must have $k+1$ vertices. Hence $w$ has exactly $k+1$
stops; we list them as
 $S^0, S^1, \ldots, S^k$. For $1\leq i \leq k$, let $E_i$ be the set of $(r-2s)$ vertices
in $F^i$ but not in any stops.
We get a partition:
$U=\left(\cup_{j=0}^iS^j\right)\cup \left( \cup_{j=1}^i E_j\right)$.
A code consists of $k$ `$($' and $k$ `$)$' is generated as follows.
When we scan the walk from left to right, if an edge appears the first time,
we append the code by a  `$($'; otherwise, we append the code by a `$)$'.
The code is  a valid sequence of $k$ pairs of parentheses. (In this case, the
number of `$\ast$'s is zero.)  It suffices to recover a walk from a
partition of $[m_k]$ and a sequence of valid parentheses.

Given a partition of $U$
$$\left(\cup_{j=0}^iS^j\right)\cup \left( \cup_{j=1}^i E_j\right)$$
and a sequence of $k$ pairs valid parentheses,
 we first  build a rooted tree $T$ as follows. At each time,
we maintain a tree $T$, a current stop $S$, a set of unused stops $\S$.
Initially $T$ contains nothing but the root stop $S_0$, $S:=S_0$,
and $\S=\{S_1,S_2,\ldots, S_k\}$.
At each time, read a symbol from the sequence.
If the symbol is  an open parenthesis, then find an $S_i$ in $\S$ with
index $i$ as small as possible, delete $S_i$ from $\S$, attach $S_i$ to $T$
as a child stop of $S$, and let $S:=S_i$;
if the symbol is ``)'',  then let $S$ point to the the parent stop of the current $S$.
Repeat this process until all symbols from the sequence are processed.

Since every closed parenthesis has a matching open parenthesis, this process
never get stuck.
When the process ends, a rooted tree $T$ on the vertex set $\{S_0,\ldots, S_{k}\}$ is
created. For $1\leq i \leq k$, let $F_i$ be the union of $E_i$ and
 two ends of $i$-th edge, which created in the process.
For example, for $k=3$, if the sequence is $(())()$, then
the corresponding good closed walk is
$$S_1F_1 S_2F_2 S_3F_2 S_2F_1 S_1F_3S_4F_3 S_1$$
where $F_1=S_1\cup S_2 \cup E_1$, $F_2=S_2\cup S_3 \cup E_2$,
and $F_3= S_4\cup S_1 \cup E_3$.

Thus, this is a bijection from $\G_k^{m_k}$ to all triples $\{U, \P, C\}$.
The number of ways to choose $m_k$ vertices is ${n\choose m_k}$.
 The number of ways to choose
these sets $S_0, S_1, \ldots, S_k, E_1,\ldots, E_k$ as  a partition of $U$  is
$${m_k \choose s, \ldots, s, r-2s,\ldots, r-2s}=\frac{m_k!}{(s!)^{k+1}((r-2s)!)^{k}}.$$
The number of   sequences of $k$ pairs valid parentheses is the
Catalan number $\frac{1}{k+1}{2k \choose k}$.  By taking product of
these three numbers, we get equation (\ref{eq:gk}). \hfill $\square$

\noindent {\bf Proof of Lemma \ref{trace}:} By equations
(\ref{eq:trace1}) and (\ref{eq:gi}), we have
$$\E(\tr(C^{t}))\leq
\sum_{i=1}^{\lfloor \frac{t}{2}\rfloor}|\G_i| p^i(1-p)^i\leq (1+o(1))
\sum_{i=1}^{\lfloor \frac{t}{2}\rfloor}a_i.$$
Here $a_i:={t-2\choose t-2i} i^{t-2i}\frac{1}{i+1}{2i\choose i} {r-s\choose s}^{t-i}
\frac{n^{m_i}p^i(1-p)^i}{(s!)^{i+1}((r-2s)!)^{i}} $. We get
\begin{eqnarray*}
  \frac{a_i}{a_{i+1}} &=& 
\frac{{t-2\choose t-2i} i^{t-2i}\frac{1}{i+1}{2i\choose i} {r-s\choose s}^{t-i}
\frac{n^{m_i}p^i(1-p)^i}{(s!)^{i+1}((r-2s)!)^{i}}
}
{{t-2\choose t-2i-2} i^{t-2i-2}\frac{1}{i+2}{2i+2\choose i+1} {r-s\choose s}^{t-i-1}
\frac{n^{m_{i+1}}p^{i+1}(1-p)^{i+1}}{(s!)^{i+2}((r-2s)!)^{i+1}}}  \\
&=&\frac{i^3(2i-1)(i+2)}{(2i+1)(t-2i)(t-2i-1)}\frac{(r-s)!}{n^{r-s}p(1-p)}\\
&<& \frac{3i^4(r-s)!}{n^{r-s}p(1-p)}.
\end{eqnarray*}

When $n^{r-s}p(1-p)\gg t^4$,  we have $a_i=o(a_{i+1})$. Thus,
$$\E(\tr(C^{t})) \leq (1+o(1)) a_{\lfloor \frac{t}{2}\rfloor}.$$
When $t=2k$, we get
$$ \E\left(\tr(C^{2k})\right) \leq   (1+o(1))
\frac{n^{s+k(r-s)}{r-s \choose s}^k}{(k+1) (s!)^{k+1}((r-2s)!)^k} \binom{2k}{k} p^k(1-p)^k.$$

For $t=2k+1$, we have
\begin{eqnarray*}
 \E(\tr(C^{2k+1}))&\leq& (1+o(1))a_k\\
&\leq& (1+o(1))\frac{2k^2n^{s+k(r-s)}{r-s \choose s}^{k+1}}{(k+1) (s!)^{k+1}((r-2s)!)^k} \binom{2k}{k} p^k(1-p)^k.
\end{eqnarray*}

Now we assume $k=o(\log (n^{r-s}p(1-p)))$. For $t=2k$,  let
$$b_k:=|\G_k^{m_k}| p^k(1-p)^k={n\choose m_k}\frac{m_k!}{(k+1)(s!)^{k+1}((r-2s)!)^{k}}\binom{2k}{k} p^k(1-p)^k.$$
It is clear that  $\E(\tr(C^{2k}))\geq b_k$. We also have
\begin{eqnarray*}
  \E(\tr(C^{2k}))-b_k&\leq& \sum_{i=1}^{k-1}|\G_i| p^i(1-p)^i +
\sum_{j=r}^{m_k-1}|\G_k^j| p^k(1-p)^k\\
&=& (1+o(1))\sum_{i=1}^{k-1}a_i +
a_k \sum_{j=r}^{m_k-1} \left(\frac{C_1 k^{C_2}}{n}\right)^{m_k-j}.
\end{eqnarray*}
Note $a_k=(1+o(1)){r-s\choose s}^k b_k$ and
$a_i=O\left(a_k \left(\frac{k^4}{n^{r-s}p(1-p)}\right)^{k-i}\right)$. We conclude
$$\E(\tr(C^{2k}))-b_k = O\left(b_k  {r-s\choose s}^k 
\left(\frac{k^4}{n^{r-s}p(1-p)} +\frac{C_1 k^{C_2}}{n}\right)
\right)=o(b_k).$$
Here we use the fact   ${r-s\choose s}^k
\left(\frac{k^4}{n^{r-s}p(1-p)}+\frac{C_1 k^{C_2}}{n}\right)=o(1)$ since $k=o(\log (n^{r-s}p(1-p)))$.
\hfill $\square$



\section{The semicircle law}
Let us review the definition of the Semicircle Law.
Let $F(x)$ be the continuous distribution function with
density $f(x)$ such that $f(x)=(2/\pi) \sqrt{1-x^2}$ when $|x| \leq
1$ and $f(x)=0$ when $|x|>1.$
Let $A$ be a Hermitian matrix of dimension $N\times N$.
The {\em empirical distribution} of the eigenvalues of $A$ is
$$F(A,x):=\frac{1}{N}|\{\mbox{ eigenvalues of } A
\mbox{ less than } x \}|.$$

We say, the empirical distribution of the eigenvalues of $A$ asymptotically
follows the Semicircle Law centered at $c$  with radius $R$ if
$F(\frac{1}{R}(A-cI), x)$ tends to $F(x)$ in probability as $N$ goes to infinity.
(In this case, we write $F(\frac{1}{R}(A-cI), x) \stackrel{p}{\to} F(x)$.)
If $c$ is the  center of the Semicircle Law, then any $c'=c+o(R)$ is also the center
of the Semicircle Law.

\begin{theorem} \label{t3}
If  $n^{r-s}p(1-p)\to \infty$,  then the empirical distribution of
the eigenvalues of $W-\E(W)$ follows the
semicircle law centered at $0$ with radius $2\sqrt{\binom{r-s}{s}
\binom{n-s}{r-s}p(1-p)}$.
\end{theorem}

{\bf Proof:}  Let $R:=2\sqrt{\binom{r-s}{s}
\binom{n-s}{r-s}p(1-p)}$,
 $C:=W-\E(W)$, and  $C_{nor}:=\frac{1}{R}C$.

To prove the theorem, we need  to show that for any
fixed $t$, the $t$-th moment of $F(C_{nor},x)$ (with $n$ goes to infinity)
is asymptotically equal to the $t$-th moment of $F(x)$. We know the
$t$-th moment of $F(C_{nor},x)$ equals $\binom{n}{s}^{-1}
\E(\tr(C_{nor}^t))$. For even $t=2k$,  the $t$-th moment of $F(x)$
is $(2k)!/2^{2k}k!(k+1)!$. For odd $t$, the $t$-th moment of $F(x)$
is 0.

In order to prove the theorem, we need to show for any fixed $k$,
$$
\frac{1}{\binom{n}{s}}
\E(\tr(C_{nor}^{2k}))=(1+o(1))\frac{(2k)!}{2^{2k}k!(k+1)!}
$$
and
$$
\frac{1}{\binom{n}{s}} \E(\tr(C_{nor}^{2k+1}))=o(1).
$$
 We know
$$
 \E(\tr(C_{nor}^{t}))=\frac{1}{R^{t}} \E(\tr(C^{t}))
$$
for any $t$.  By Lemma \ref{trace}, we have
$$
\E(\tr(C^{2k})) = (1+o(1)) \frac{n^{s+k(r-s)}}{(k+1) (s!)^{k+1} ((r-2s)!)^k} \binom{2k}{k} p^k(1-p)^k.
$$
Then
$$
\frac{1}{\binom{n}{s}} \E(\tr(C_{nor}^{2k}))=(1+o(1))\frac{(2k)!}{2^{2k}k!(k+1)!}
$$
as desired.

By Lemma \ref{trace} again, we have
$$
\E(\tr(C^{2k+1}))=O\left (\frac{2k^2n^{s+k(r-s)}p^k(1-p)^k{r-s\choose s}^{k+1}}{(k+1)(s!)^{k+1}( (r-2s)!)^k}{2k\choose k}\right).
$$
Thus
\[\frac{1}{\binom{n}{s}} \E(\tr(C_{nor}^{2k+1}))
= O\left(\frac{2k^2{2k\choose k} {r-s\choose s}^{k+1}}{2^{2k}(k+1)R}\right)\\
=o(1).\]
Here $k$ is any constant but $R\to\infty$.
 The theorem is proved. \hfill $\square$

The following Lemma is useful to derive  the Semicircle Law from one
matrix to the other.
 \begin{lemma} \label{semicircle}
   Let $A$ and $B$ be two $(N\times N)$-Hermitian matrices. Suppose that
the empirical distribution of the eigenvalues of $A$
follows the Semicircle Law centered at $c$
with radius $R$. If either $\|B\|=o(R)$ or the rank of $B$ is $o(N)$, then
the empirical distribution of the eigenvalues of $A+B$
also follows the Semicircle Law centered at $c$
with radius $R$.
 \end{lemma}
{\bf Proof:}
It suffices to show $F(\frac{1}{R}(A+B-cI),x) \stackrel{p}{\to} F(x)$.
First we assume $\|B\|=o(R)$. By Lemma \ref{difference}, for $1\leq k\leq N$,
we have
$$\left |\mu_k\left(\frac{1}{R}(A+B-cI)\right)- \mu_k\left(\frac{1}{R}(A-cI)\right)\right |
\leq \frac{\|B\|}{R}=o(1).$$
Hence
$$F \left(\frac{1}{R}(A-cI),x-\frac{\|B\|}{R}\right)\leq
F\left(\frac{1}{R}(A+B-cI),x\right)
\leq F\left(\frac{1}{R}(A-cI),x+\frac{\|B\|}{R}\right).$$
Since $\|B\|=o(R)$, we have $F \left(\frac{1}{R}(A-cI),x-\frac{\|B\|}{R}\right) \stackrel{p}{\to} F(x)$ and $F \left(\frac{1}{R}(A-cI),x+\frac{\|B\|}{R}\right) \stackrel{p}{\to} F(x)$.
By the Squeeze theorem,  we have $F(\frac{1}{R}(A+B-cI),x) \stackrel{p}{\to} F(x)$.

Now we assume ${\rm rank}(B)=o(N)$.
Let $U$ be the kernel of $B$ (i.e. $B|_U=0$); $U$ has co-dimension ${\rm rank(B)}$.
Let $Z:=\frac{1}{R}(A-cI)|_U=\frac{1}{R}(A+B-cI)|_U$.
By Cauchy's interlace theorem \cite{cauchy},
for $1\leq j\leq N-{\rm rank}(B)$,
we have
\begin{eqnarray*}
  \mu_j\left(\frac{1}{R}(A-cI)\right) &\leq& \mu_j(Z)
\leq  \mu_{j+{\rm rank}(B)}\left(\frac{1}{R}(A-cI)\right),\\
  \mu_j\left(\frac{1}{R}(A+B-cI)\right) &\leq&
\mu_j(Z) \leq  \mu_{j+{\rm rank}(B)}\left (\frac{1}{R}(A+B-cI)\right).
\end{eqnarray*}
Thus,  for ${\rm rank}(B) +1\leq j \leq N-{\rm rank}(B)$, we have
$$\mu_{j-{\rm rank}(B)}\left(\frac{1}{R}(A-cI)\right)\leq
\mu_{j}\left (\frac{1}{R}(A+B-cI)\right)
\leq \mu_{j+{\rm rank}(B)}\left(\frac{1}{R}(A-cI)\right).$$
It implies
$$F\left(\frac{1}{R}(A-cI),x\right) -
\frac{{\rm rank}(B)}{N} \leq
F\left(\frac{1}{R}(A+B-cI),x\right)\leq F\left(\frac{1}{R}(A-cI),x\right) +
\frac{{\rm rank}(B)}{N}.$$

Since ${\rm rank}(B)=o(N)$, we have
$F\left(\frac{1}{R}(A-cI),x\right) \pm\frac{{\rm rank}(B)}{N} \stackrel{p}{\to} F(x)$.
By the Squeeze theorem,  we have $F(\frac{1}{R}(A+B-cI),x) \stackrel{p}{\to} F(x)$.
\hfill $\square$

{\bf Proof of Theorem \ref{t2}}:
Recall $$\L^{(s)}(K^r_n)-\L^{(s)}(H^r(n,p))=M_1+M_2+M_2+M_4.$$
We can write $\L^{(s)}(H^r(n,p))$ as $-M_2 + \left(1-\frac{(-1)^s}{{n\choose s}}\right) I + B_1-M_3-M_4-M_1 $,
where $B_1=\L^{(s)}(K^r_n)-\left(1-\frac{(-1)^s}{{n\choose s}}\right) I$.

By Theorem \ref{t3}, the empirical distribution of the spectrum of $W-\E(W)$
follows the Semicircle Law centered at $0$
 with radius $(2+o(1))\sqrt{\binom{r-s}{s} \binom{n-s}{r-s}p(1-p)}$.
Since $M_2=\frac{1}{{r-s\choose s}d}(W-\E(W))$, $ \left(1-\frac{(-1)^s}{{n\choose s}}\right) I -M_2$ follows
the Semicircle Law centered at $c:=1-\frac{(-1)^s}{{n\choose s}}$ with radius
$R:=(2+o(1))\sqrt{\frac{1-p}{\binom{r-s}{s} \binom{n-s}{r-s}p}}$. Note
$\frac{(-1)^s}{{n\choose s}}=o(R)$. We can change the center to $1$.

By Theorem \ref{t0}, $\L^{(s)}(K^r_n)$ has an eigenvalue
$1-(-1)^s\frac{{n-s\choose s}}{{n\choose s}}$ with multiplicity
${n\choose s}-{n\choose s-1}$. Thus $B_1$ has rank ${n\choose s-1}=
o\left({n\choose s}\right)$. We also observe
that $M_4$ has rank at most $2$, $\|M_1\| =O\left(\frac{\sqrt{(1-p)\log N}}{d}\right)=o(R)$,
and $\|M_3\|=O\left(\frac{\sqrt{\log N}}{n\sqrt{d}}\right)=o(R)$.
 Here we notice $d\gg \log^{1/3} n$ and $1-p\gg \frac{\log n}{n^2 d^2}$.

By Lemma \ref{semicircle}, the matrices $B_1$, $M_1$, $M_3$, and $M_4$ will not
affect the Semicircle Law. The proof of this Lemma is finished.
\hfill $\square$

\end{document}